%% file: DissBelowScaling.tex
\documentclass[11pt,headsepline,abstracton,titlepage,oneside]{scrartcl}
\usepackage[a4paper]{geometry}
\usepackage{amsmath,amsthm,amsfonts,amssymb}
\usepackage{graphicx,color}
\usepackage{makeidx}
\newtheorem{thm}{Theorem}[section]
\newtheorem{lem}[thm]{Lemma}
\newtheorem{prop}[thm]{Proposition}
\newtheorem{cor}[thm]{Corollary}

\newtheorem{df}{Definition}[section]
\theoremstyle{remark}
\newtheorem{rem}{Remark}[section]
\newtheorem{expl}[rem]{Example}
\numberwithin{equation}{section}
\numberwithin{figure}{section}
\def\<{\left<}
\def\>{\right>}
\def\R{\mathbb R}
\def\C{\mathbb C}
\def\N{\mathbb N}

\def\d{\mathrm d}
\def\D{\mathrm D}

\newcommand{\bx}{\boldsymbol{\square}}
\let\into\hookrightarrow

\DeclareMathOperator{\supp}{supp}

\DeclareMathOperator*{\esssup}{ess\,sup}

\DeclareMathOperator*{\slim}{s-lim}

\DeclareMathOperator{\dist}{dist}

\newcounter{prgrph}[section]

\newcommand{\eIndex}[1]{#1\index{#1}}
\begin{document}
\titlehead{TU Bergakademie Freiberg\\
  Fakult\"at f\"ur Mathematik und Informatik\\
  Institut f\"ur Angewandte Analysis\\
  09596 Freiberg, GERMANY}
\author{Jens Wirth\thanks{\tt mailto:wirth@math.tu-freiberg.de}}
\title{Asymptotic properties of solutions to weakly dissipative
wave equations below scaling}
\def\today{July 27, 2004}
\begin{titlepage}
\maketitle
\end{titlepage}

%%%%%%%%%%%%%%%%%%%%%%%%%%%%%%%%%%%%%%%%%%%%%%%%%
%%
%%    SEC: INTRO
%%
%%%%%%%%%%%%%%%%%%%%%%%%%%%%%%%%%%%%%%%%%%%%%%%%%

\setcounter{section}{-1}
\section{Introduction}

Cauchy problems for weakly dissipative wave equations
\begin{equation}\label{chap3:eq:CP}
  u_{tt}-\Delta u+b(t)u_t=0, \qquad u(0,\cdot)=u_1,\quad \D_tu(0,\cdot)=u_2,
\end{equation}
where the coefficient $b=b(t)$ is assumed to be positive and tends to zero as
$t$ tends to infinity, provide an important model problem for the study of 
asymptotic behaviours and the influence of lower order terms on them. We refer to
\cite{WirOverview} for an exposition of results in that case or to
the papers of A.~Matsumura \cite{Mat77}, H.~Uesaka \cite{Ues79} or
F.~Hirosawa and H.~Nakazawa \cite{HN03} for known energy decay estimates. 

In special cases solution representations are explicitly known and thus also the question
of sharpness of estimates can be answered. For properties of free waves we refer to
the book of R.~Racke \cite{Racke} or the classical paper of R.S.~Strichartz,
\cite{Str69}. On the other hand the case of scale invariant weak
dissipation, $b(t)\sim t^{-1}$ as $t\to\infty$, is closely related to Bessel's
differential equation and considered by the author in \cite{Wir02}. 

The purpose of this paper is to prove structural properties of representations of solutions
in the case below scaling, i.e. we assume that the coefficient function
decays faster than $t^{-1}$. This case is closely related to free waves. In case of integrable
coefficients it is known from the works of K.~Mochizuki and coauthors \cite{Moc76}, \cite{MN96} 
or from the preprint \cite{Wir02b} that the solutions are asymptotically free and a scattering theory in the sense
P.D.~Lax and R.S.~Philipps \cite{LP73} can be obtained.

We will employ the translation invariance of the Cauchy problem. This 
implies that a partial Fourier transform with respect to the spatial
variables may be used to reduce the partial differential equation
in $u(t,x)$ to an ordinary differential equation for $\hat u(t,\xi)$ parameterised
by the frequency parameter $|\xi|$,
$$
  \hat u_{tt}+b(t)\hat u_t+|\xi|^2\hat u=0.
$$
Its solution can be represented in the form
$$
  \hat u(t,\xi)=\Phi_1(t,\xi)\hat u_1+\Phi_2(t,\xi)\hat u_2
$$
in terms of the Cauchy data $u_1$ and $u_2$ with suitable
functions (Fourier multiplier) $\Phi_1$ and $\Phi_2$. 

In general estimates for $\Phi$ are complicated to obtain, so the natural
starting point is to rewrite the second order equation as system for the 
energy $(|\xi|\hat u,\D_t u)^T$ or a related vector and to use
a diagonalization technique to simplify the structure and to estimate
its fundamental solution. 

Main results of this paper are the solution representation of Theorem 
\ref{thm:II:Zhyp_highreg_solrep} together with its consequences for the
$L^p$--$L^q$ decay, Theorem \ref{thm:II:LpLq-decay}.
Furthermore the sharpness of these results follows from a modified scattering
theory given by Theorem \ref{thm:II:mod_scattering}.

\section{Representation of Solutions}

We start by giving examples. The technique used in this paper allows to consider the 
following coefficient functions.
\begin{expl}
  Scale invariant dissipation, $b(t)=\frac\mu{1+t}$ is contained for $\mu\in[0,1)$.
  For the general case $\mu>0$ we refer to \cite{Wir02}. The exceptional value $\mu=1$ 
  is critical for estimates of the solution itself and related to an exceptional
  behaviour of Bessel functions.
\end{expl}
\begin{expl} Let $\mu>0$ and $n\ge 1$. Then we can consider 
   $$ b(t)=\frac{\mu}{(1+t)\log(e+t)\cdots\log^{[n]}(e^{[n]}+t)}. $$
  Here we denote 
  by $\log^{[n]}$ the  $n$ times iterated logarithm, $\log^{[0]}(\tau)=\tau$ 
  and $\log^{[k+1]}(\tau)=\log\log^{[k]}(\tau)$, and by $e^{[n]}$ similar the iterated exponential
  $e^{[0]}=1$ and $e^{[k+1]}=e^{(e^{[k]})}$.
\end{expl}
\begin{expl}
 Monotonicity of the coefficient is not essential. We can also consider coefficient functions
 like
 $$
   b(t)=\frac{2+\cos(\alpha\log(e+t))}{4(e+t)}
 $$
 with $\alpha\in\R$.
\end{expl}

\paragraph{Assumptions.}
We make the following assumptions on the coefficient function $b=b(t)$
\begin{description}
\item[(A1)]\index{Assumption!(A1)}
    positivity: $b(t)\ge 0$, 
\item[(A2)]\index{Assumption!(A2)}symbol-like estimates:
   for all $k\in\N$ it holds  
  $$ \left|\frac{\d^k}{\d t^k} b(t)\right|\le C_k \left(\frac1{1+t}\right)^{1+k},  $$ 
\item[(A3)]\index{Assumption!(A3)}
    it holds $\limsup_{t\to\infty} tb(t)<1$.
\end{description}
The last assumption is necessary to exclude the above mentioned critical case
$b(t)=\frac1{1+t}$.

\paragraph{Basic ideas, zones.} If we consider the vector valued function 
$\tilde U=(|\xi|\hat u,\D_t\hat u)^T$, the energy, we obtain the system
\begin{equation}
  \D_t \tilde U=\begin{pmatrix}&|\xi|\\|\xi|&ib(t) \end{pmatrix} \tilde U.
\end{equation}
In order to understand the properties of its solutions we relate the size of 
$b(t)$ with the size of $|\xi|$. This leads to the
introduction of different zones in the extended phase space.

Therefore, we will use the implicitly defined function
\begin{equation}\label{eq:zones:high-reg}
  (1+t_\xi)|\xi|=N
\end{equation}
with suitable constant $N$ and define the
\begin{equation}
  Z_{hyp}(N)=\{(t,\xi)|t\ge t_\xi\},\qquad\qquad Z_{diss}(N)=\{(t,\xi)|0\le t\le t_\xi\}.
\end{equation} 
Inside $Z_{diss}$ the positivity of the dissipation $b(t)$ is essential and
we will use a reformulation as integral equation to solve the system.

In $Z_{hyp}$ the non-diagonal entries are dominating and we use a diagonalization scheme
to transform the system to a diagonally dominated one.
In order to do this we employ the symbol estimate of Assumption (A2) 
and introduce symbol classes in the hyperbolic zone $Z_{hyp}$. This is related to the approach
in \cite{RY00}. 

\begin{df}\label{def:Symb_hyp}
  The time-dependent Fourier multiplier $a(t,\xi)$ belongs to the \eIndex{hyperbolic symbol class}
  $S_N\{m_1,m_2\}$ if it is supported in the hyperbolic zone and satisfies the symbol estimates
  \begin{equation}\label{eq:symb_est_hyp}
    \left|\D_t^k\D_\xi^\alpha a(t,\xi)\right|\le 
    C_{k,\alpha} |\xi|^{m_1-|\alpha|}\left(\frac1{1+t}\right)^{m_2+k}
  \end{equation}
  for all $(t,\xi)\in Z_{hyp}(N)$ and all $k\in\N$, $\alpha\in\N^n$.
\end{df}

Thus using \eqref{eq:zones:high-reg} we can conclude the essential embedding rule
\begin{equation}
   S_{N}\{m_1-k,m_2+\ell\}\into S_N\{m_1,m_2\}\qquad \forall\; \ell\ge k\ge0.
\end{equation}
Definition \ref{def:Symb_hyp} extents immediately to matrix-valued Fourier multiplier. The rules
of the symbolic calculus are obvious and collected in the following proposition.
\begin{prop}\label{prop:II:calc_rules}
\begin{enumerate}
\item $S_N\{m_1,m_2\}$ is a vector space.
\item $S_N\{m_1,m_2\}\cdot S_N\{m_1',m_2'\}\into S_N\{m_1+m_1',m_2+m_2'\}$
\item $\D_t\D_\xi^\alpha S_N\{m_1,m_2\}\into S_N\{m_1-|\alpha|,m_2+k\}$
\item $S_N\{-1,2\}\into L^\infty_\xi L^1_t$
\end{enumerate}
\end{prop}

For later use we define the larger symbol classes $S_N^{\ell_1,\ell_2}\{m_1,m_2\}$ of
symbols with restricted smoothness. For these symbols we assume the symbol estimate
\eqref{eq:symb_est_hyp} for $k=0,1,\ldots\ell_1$ and $|\alpha|\le\ell_2$. These symbols 
with restricted smoothness can be used to deduce mapping properties in $L^p$ spaces\footnote{We are 
speaking about Fourier multiplier only, so no essential difficulties can arise by this 
lack of smoothness.}. We give one auxiliary result following directly from 
Marcinkiewicz multiplier theorem, \cite[Chapter IV.3, Theorem 3]{Ste70}. 

\begin{prop} 
  Each $a\in S^{0,\lceil\frac n2\rceil}_N\{0,m\}$ gives rise to an operator
  $a(t,\D):L^p\to L^p$ for all $p\in(1,\infty)$ with norm estimate
  $$ ||a(t,\D)||_{p\to p}\lesssim \left(\frac1{1+t}\right)^m. $$
\end{prop}

\paragraph{Formulation in system form.}
We do not use the vector function $\tilde U$ mentioned before. The main reason is
that we prefer a representation, where we can extract also estimates for the 
solution itself. Furthermore for small frequencies we need information on the size
of the solution in order to obtain sharp decay results.

We consider the \eIndex{micro-energy}
\begin{equation}\label{eq:U_def}
  U=(h(t,\xi)\hat u,\D_t\hat u)^T
\end{equation}
with
\begin{equation}
  h(t,\xi)=\frac N{1+t}\phi_{diss,N}(t,\xi)+|\xi|\phi_{hyp,N}(t,\xi).
\end{equation}
Here and thereafter we denote by $\phi_{diss,N}(t,\xi)$ the characteristic function of the
dissipative zone and by $\phi_{hyp,N}(t,\xi)$ the characteristic function of the hyperbolic zone.  

This micro-energy is related to the energy operator $\mathbb E(t,\D)$,
\begin{equation}
 \mathbb E(t,\D):(\<\D\>u_1,u_2)^T\mapsto (|\D|u(t,\cdot),\D_t u(t,\cdot))^T,
\end{equation}
also considered in \cite{Wir02}. For all $t$ we estimate some
kind of $H^1$-norm, but uniformly for all $t$ the estimate can only be seen as $\dot H^1$ estimate.

The aim is to prove estimates and structural properties for the fundamental solution
$\mathcal E(t,s,\xi)$ to the corresponding system $\D_tU=A(t,\xi)U$.

\subsection{The dissipative zone}
The essential idea is to write the problem as an integral equation. 
In the dissipative zone the micro-energy \eqref{eq:U_def} simplifies to
$$
  U=\left(\frac N{1+t}\hat u,\D_t\hat u\right)^T
$$
and thus we have to solve the system
\begin{equation}\label{eq:II:Zdiss_CP}
  \D_t\mathcal E(t,s,\xi)=A(t,\xi)\mathcal E(t,s,\xi)
  =\begin{pmatrix}\frac i{1+t} & \frac N{1+t} \\ \frac{(1+t)|\xi|^2}N & ib(t)\end{pmatrix}
  \mathcal E(t,s,\xi),\qquad \mathcal E(s,s,\xi)=I.
\end{equation}

We rewrite it as an integral equation. Let therefore $\Lambda(t,s)$ be the
fundamental solution of the diagonal part of this system,
\begin{equation}
  \Lambda(t,s)=\begin{pmatrix}\frac{1+s}{1+t}&\\&\frac{\lambda^2(s)}{\lambda^2(t)}\end{pmatrix},
\end{equation}
$\lambda(t)=\exp\left\{\frac12\int_0^tb(\tau)\d\tau\right\}$,
such that the application of Duhamel's principle yields
\begin{equation}
  \mathcal E(t,s,\xi)=\Lambda(t,s)
  +i\int_s^t\Lambda(t,\tau)R(\tau,\xi)\mathcal E(\tau,s,\xi)\d\tau
\end{equation}
with $R(t,\xi)=A(t,\xi)-\mathrm{diag}\;A(t,\xi)$. Assumption (A1) guarantees that the
entries of $\Lambda(t,s)$ are monotonous functions.

We claim that the dominating entry of $\Lambda(t,s)$ describes the behaviour of 
$\mathcal E(t,s,\xi)$. If we consider $\Lambda(s,t)\mathcal E(t,s,\xi)$ we get the
integral equation
\begin{equation}
  \Lambda(s,t)\mathcal E(t,s,\xi)=I
  +i\int_s^t\underbrace{\Lambda(s,\tau)R(\tau,\xi)\Lambda(\tau,s)}_{K(\tau,s,\xi)}
     \big(\Lambda(s,\tau)\mathcal E(\tau,s,\xi)\big)\d\tau
\end{equation}
with kernel $K(t,s,\xi)$. We will prove well-posedness of this equation in $L^\infty$
over the dissipative zone using Theorem \ref{thm:B4}. It holds
\begin{equation}
  K(\tau,s,\xi)=\begin{pmatrix}&\frac{N\lambda^2(s)}{\lambda^2(\tau)(1+s)} \\
   \frac{\lambda^2(\tau)(1+s)}{\lambda^2(s)N}|\xi|^2\end{pmatrix}.
\end{equation}

Therefore we have to relate the auxiliary function $\lambda^2(t)$
to $t$. This is done using Assumption (A3) and the following consequence of it.

\begin{prop}\label{prop:8ab}
 Assumptions (A1), (A3) imply 
  $$\int_0^t\frac{\d\tau}{\lambda^2(\tau)} \sim \frac t{\lambda^2(t)}$$ 
  and $\frac t{\lambda^2(t)}$ is monotonous increasing for large $t$.
\end{prop}
\begin{proof}
 Integration by parts yields
$$ \int_0^t \frac{\d\tau}{\lambda^2(\tau)}=\frac t{\lambda^2(t)}+\int_0^t \frac{\tau b(\tau)}{\lambda^2(\tau)}\d\tau. $$
On the one hand the right hand side is larger than $t\lambda^{-2}(t)$ by Assumption (A1). On  the other hand we conclude
from $tb(t)\leq c<1$ for $t>t_0$ that
$$ \int_0^t \frac{\tau b(\tau)}{\lambda^2(\tau)}\d\tau
\leq\int_0^{t_0} \frac{\tau b(\tau)}{\lambda^2(\tau)}\d\tau+c\int_{t_0}^t \frac{\d\tau}{\lambda^2(\tau)}
\leq C+c\int_{0}^t \frac{\d\tau}{\lambda^2(\tau)}$$
and the statement follows
$$ \int_0^t \frac{\d\tau}{\lambda^2(\tau)}\leq \frac1{1-c} \left(C+\frac t{\lambda^2(t)}\right)\lesssim \frac{t}{\lambda^2(t)}.$$
For small $t$ the statement can be concluded from $\lambda^2(t)\sim 1$.

Monotonicity is a consequence of
$$ \frac{\d}{\d t}\frac t{\lambda^2(t)}=\frac{1-tb(t)}{\lambda^2(t)}$$
and $tb(t)<1$ for $t\gg1$.
\end{proof}

Especially under Condition (A3) we have $\lambda^2(t)\lesssim (1+t)$.

\begin{lem}\label{lem:Zdiss1}
  Assume (A1) and (A3). Then 
  \begin{equation}\label{eq:E_bound}
  ||\mathcal E(t,s,\xi)||\lesssim \frac{\lambda^2(s)}{\lambda^2(t)},\qquad t_\xi\ge t\ge s.
  \end{equation}
\end{lem}
\begin{proof}
  It remains to check the conditions for Theorem \ref{thm:B4}. It holds
  $$
  \frac{\lambda^2(s)}{1+s}\int_s^{t_\xi} \frac{\d\tau}{\lambda^2(\tau)}
    \sim \frac{\lambda^2(s)(1+t_\xi)}{\lambda^2(t_\xi)(1+s)}- 1
    \lesssim \frac{\lambda(s)}{\lambda(t_\xi)}\frac N{(1+s)|\xi|}\lesssim 1
  $$
  using the definition of the zone together with
  $$
   |\xi|^2\frac{1+s}{\lambda^2(s)} \int_s^t\lambda^2(\tau)\d\tau
   \lesssim  |\xi|^2\frac{1+s}{\lambda^2(s)}\big(t_\xi\lambda^2(t_\xi)-s\lambda^2(s)\big)
   \lesssim 1
  $$
  following from the monotonicity of $\frac{t}{\lambda^2(t)}$ for large $t$.
\end{proof}

\paragraph{Further results for higher order derivatives.} In order to perform a perfect diagonalization
in the hyperbolic zone it is essential to find symbol estimates for $\mathcal E(t_\xi,0,\xi)$ for 
$|\xi|\le N$. 
\begin{lem}
Assume that (A1) and (A3) holds.
Then for $|\xi|\leq N$ and all $\alpha\in\mathbb N^n$ the symbol like estimate
  $$ ||\D_\xi^\alpha \mathcal E(t_\xi,0,\xi)||\leq C_\alpha \frac1{\lambda^ 2(t_\xi)}|\xi|^{-\alpha} $$
is valid.
\end{lem}

\begin{proof}
  It holds $\D_t\mathcal E=A\mathcal E$ with 
  $$ A(t,\xi)=\begin{pmatrix}&\frac N{1+t}\\\frac{1+t}N|\xi|^2&ib(t)\end{pmatrix},\qquad 
     ||A(t,\xi)||\lesssim \frac1{1+t}. $$
  Thus for $|\alpha|=1$ we get
  $$ \D_t\D_\xi^\alpha \mathcal E=\D_\xi^\alpha (A\mathcal E)
     =(\D_\xi^\alpha A)\mathcal E+A(\D_\xi^\alpha \mathcal E) $$
  or using Duhamel's formula
  together with the initial value 
  $\D_\xi^\alpha\mathcal E(0,0,\xi)=0$
  $$ \D_\xi^\alpha \mathcal E(t,0,\xi)= 
    \int_0^t \mathcal E(t,\tau,\xi) (\D_\xi^\alpha A(\tau,\xi))\mathcal E(\tau,0,\xi)\d\tau. $$
  Now the known estimates for $\mathcal E(t,s,\xi)$, \eqref{eq:E_bound},
  imply together with $||\D_\xi^\alpha A(t,\xi)||\lesssim 1$ the desired statement
  $||\D_\xi^\alpha\mathcal E(t,0,\xi)||\lesssim t\lesssim |\xi|^{-1}$ for $(t,\xi)\in Z_{diss}(N)$ .
  
  For $|\alpha|=\ell>1$ we use Leibniz rule to get similar representations 
  containing all derivatives of order 
  less than $|\alpha|$ under the integral and use induction over $\ell$.

  The estimates 
  $$ (\D_\xi^{\alpha_1}A)(\D_\xi^{\alpha_2}\mathcal E)
     \lesssim \frac1{\lambda^2(t)}|\xi|^{1-|\alpha_1|-|\alpha_2|} $$
  for $|\alpha_1|+|\alpha_2|\leq \ell$, formula \eqref{eq:E_bound} and the first 
  statement we conclude
  $$ ||\D_\xi^\alpha \mathcal E(t,\xi)||\lesssim \int_0^t   
          \frac1{\lambda^2(\tau)}|\xi|^{1-\ell}\d\tau
    \lesssim \frac1{\lambda^2(t)} |\xi|^{-|\alpha|} $$
  with Proposition \ref{prop:8ab}.1.
  By application of the equation itself we get estimates for time derivatives. This means
  $$ ||\D_t^k\D_\xi^\alpha \mathcal E(t,\xi)||\lesssim 
     \frac1{\lambda^2(t)} \left(\frac1{1+t}\right)^k |\xi|^{-|\alpha|} $$
  and together with
  \begin{equation}\label{eq:txi_est}
    |\D_\xi^\alpha t_\xi|\lesssim |\xi|^{-1-|\alpha|},\qquad |\xi|\le N
  \end{equation}
  the statement follows.
\end{proof}

This result can be reformulated as the symbol estimate
$$
   \lambda^2(t_\xi)\mathcal E(t_\xi,0,\xi)\in\dot S^0,
$$
where 
$$ \dot S^k=\{\;m\in C^\infty(\R^n\setminus\{0\}\;|\;
    \forall\alpha\;:\;|\D^\alpha m(\xi)|\le C_\alpha |\xi|^{k-|\alpha|}\;\}$$ 
denotes the homogeneous symbol class of order $k$. Thus as
consequence of Marcinkiewicz multiplier theorem the Fourier multiplier
with symbol $\lambda^2(t_\xi)\mathcal E(t_\xi,0,\xi)$ maps $L^p$ into
$L^p$ for all $p\in(1,\infty)$.

\paragraph{Remark.} For all proven results on $\mathcal E(t,s,\xi)$ within the
dissipative zone we do not need any estimates for the 
derivatives of the coefficient function $b(t)$. The essential point is always the 
positivity of the coefficient and therefore the monotonicity of the auxiliary function
$\lambda(t)$.

\subsection{The hyperbolic zone}\label{sec:II:high_reg}

\paragraph{Diagonalization.}
 We use the special symbol classes introduced by Definition \ref{def:Symb_hyp}. Remark
that it holds $|\xi|\phi_{hyp,N}\in S_N\{1,0\}$ and by Assumption (A2) 
also $b(t)\phi_{hyp,N}\in S_N\{0,1\}$. For the further calculations we omit the 
characteristic function $\phi_{hyp,N}$.

Thus we consider
$$
  U=(|\xi|\hat u,\D_t\hat u)^T
$$
with
$$
  \D_t U=A(t,\xi)U=\begin{pmatrix} & |\xi| \\ |\xi| & ib(t) \end{pmatrix} U.
$$

We apply two transformations to this system. In a first step we diagonalize
the homogeneous principle part. After that we perform further diagonalization steps
to make the remainder belong to a sufficiently nice symbol class.

\paragraph{Step 1.} We denote by $M$ the matrix
\begin{equation}
  M=\begin{pmatrix}1&-1\\1&1\end{pmatrix}
\end{equation}
consisting of eigenvectors of the homogeneous principle part of $A(t,\xi)$
with inverse
\begin{equation}
  M^{-1}=\frac12\begin{pmatrix}1&1\\-1&1\end{pmatrix}.
\end{equation}
Then for $U^{(0)}=M^{-1}U$ we get the system
\begin{equation}
  \D_t U^{(0)}=\big( \mathcal D(\xi)+R(t)\big) U^{(0)}
\end{equation}
with 
\begin{equation}\label{eq:coef_step1}
  \mathcal D(\xi) = \begin{pmatrix}|\xi|&\\&-|\xi| \end{pmatrix},\qquad\qquad
  R(t) = \frac{ib(t)}2 \begin{pmatrix}1&1\\1&1\end{pmatrix}.
\end{equation}

\paragraph{Step $k+1$.} We construct recursively the diagonalizer $N_k(t,\xi)$ of order $k$.
Let
$$ N_k(t,\xi)= \sum_{j=0}^k N^{(j)}(t,\xi),\qquad F_k(t,\xi)=\sum_{j=0}^k F^{(j)}(t,\xi) $$
where $N^{(0)}=I$, $B^{(0)}=R(t)$ and $F^{(0)}=\mathrm{diag}\;B^{(0)}=F_0(t)$.
 
The construction goes along the following scheme. Remark that $F_0$ is a multiple of $I$.
Then we set
\begin{align*}
      &F^{(j)}=\mathrm{diag }\;B^{(j)},\\
      &N^{(j+1)} = \begin{pmatrix} & -B^{(j)}_{12}\, /\, 2|\xi|\\ 
        B^{(j)}_{21}\,/\,2|\xi| & \end{pmatrix}, \\
      &B^{(j+1)} = (\D_t-\mathcal D-R)N_{j+1}-N_{j+1}(\D_t-\mathcal D-F_j).
\end{align*}
Now we prove by induction that $N^{(j)}\in S_N\{-j,j\}$ and $B^{(j)}\in S_N\{-j,j+1\}$.  
For $j=0$ we know
$$
  F^{(0)}\in S_N\{0,1\},\quad N^{(1)}\in S_N\{-1,1\},\quad B^{(1)}\in S_N\{-1,2\},
$$ 
the last one follows from the representation $B^{(1)}=\D_t N^{(1)}-(R-F^{(0)})N^{(1)}$. 

For $j\ge1$ we apply the principle of induction. Assume we know $B^{(j)}\in S_N\{-j,j+1\}$. 
Then by definition of $N^{(j+1)}$ we have from $|\xi|^{-1}\in S_N\{-1,0\}$ 
that $N^{(j+1)}\in S_N\{-j-1,j+1\}$ and $F^{(j)}\in S_N\{-j,j+1\}$. 
Moreover,
\begin{align*}
     B^{(j+1)}&=(\D_t-\mathcal D-R)(\sum_{\nu=0}^{j+1} N^{(\nu)})
        -(\sum_{\nu=0}^{j+1} N^{(\nu)})(\D_t-\mathcal D-\sum_{\nu=0}^{j} F^{(j)})\\
     &=B^{(j)}+[N^{(j+1)},\mathcal D]-F^{(m)}+\D_tN^{(j+1)}+RN^{(j+1)}\\
      &\qquad\qquad\qquad\qquad\qquad\qquad+N^{(j+1)}\sum_{\nu=0}^j F^{(\nu)} 
    -(\sum_{\nu=1}^{j+1}N^{(\nu)})F^{(j)}.
\end{align*}
Now $B^{(j)}+[N^{(j+1)},\mathcal D]-F^{(j)}=0$ for all $j$. 
The sum of the remaining terms belongs to the symbol class $S_N\{-j-1,j+2\}$. Hence
$B^{(j+1)}\in S_N\{-j-1,j+2\}$.

Now the definition of $B^{(k)}$ implies the operator identity
\begin{equation}\label{eq:II:perf_diag}
  \big(\D_t-\mathcal D(\xi)-R(t)\big)N_k(t,\xi)
  =N_k(t,\xi)\big(\D_t-\mathcal D(\xi)-F_{k-1}(t,\xi)\big) \qquad \mod S_N\{-k,k+1\}.
\end{equation}

Thus we have constructed the desired diagonalizer if we can show that the matrix
$N_k(t,\xi)$ is invertible on $Z_{hyp}(N)$ with uniformly bounded inverse. But
this follows from $N_k-I\in S_N\{-1,1\}$ by the choice of a sufficiently large 
zone constant $N$. Indeed, we have
$$ ||N_k-I|| \leq C \frac1{|\xi|} b(t) \leq C' \frac1{|\xi|(1+t)}\to0
   \qquad\text{as $N\to\infty$}. $$

Thus with the notation $R_k(t,\xi)=-N_k^{-1}(t,\xi)B^{(k)}(t,\xi)$ we have proven the 
following lemma.

\begin{lem}\label{lem:II:perf_diag}
  Assume (A1) and (A2).
  
  For each $k\in\N$ there exists a zone constant $N$ and matrix valued symbols
  \begin{itemize}
  \item $N_k(t,\xi)\in S_N\{0,0\}$ invertible for all $(t,\xi)$ and
    with $N_k^{-1}(t,\xi)\in S_N\{0,0\}$
  \item $F_{k-1}(t,\xi)\in S_N\{0,1\}$ diagonal with 
    $F_{k-1}(t,\xi)-\frac{ib(t)}2I\in S_N\{-1,2\}$
  \item $R_k(t,\xi)\in S_N\{-k,k+1\}$, 
  \end{itemize}
  such that the (operator) identity
  \begin{equation}
  \big(\D_t-\mathcal D(\xi)-R(t)\big)N_k(t,\xi)
  =N_k(t,\xi)\big(\D_t-\mathcal D(\xi)-F_{k-1}(t,\xi)-R_k(t,\xi)\big)
  \end{equation}
  holds for all $(t,\xi)\in Z_{hyp}(N)$.
\end{lem}

\paragraph{Restricted smoothness assumptions.}
It is not necessary to know this lemma for all $k\in\N$, thus we can replace Assumption (A2)
by 
\begin{description}
\item[(A2)$_\ell$]\index{Assumption!(A2)}
  for all $k=0,1,\ldots,\ell$ it holds  
  $$ \left|\frac{\d^k}{\d t^k} b(t)\right|\le C_k \left(\frac1{1+t}\right)^{1+k}, $$ 
\end{description}
i.e. the corresponding one with finite smoothness. This gives 
$b(t)\phi_{hyp,N}\in S_N^{\ell,\infty}\{0,1\}$ and the above lemma remains true under 
Assumption (A2)$_\ell$ for all $k\le\ell$ and with matrices
$N_k\in S_N^{\ell-k+1,\infty}\{0,0\}$, $F_{k-1}\in S_N^{\ell-k+1,\infty}\{0,1\}$
and $R_k\in S_N^{\ell-k,\infty}\{-k,k+1\}$.  

\paragraph{Remarks on perfect diagonalization.} Lemma \ref{lem:II:perf_diag} can be 
understood as perfect diagonalization of the original system. If we define $F(t,\xi)$
as asymptotic sum of the $F^{(k)}(t,\xi)$, 
\begin{equation}
  F(t,\xi)\sim \sum_{k=0}^\infty F^{(k)}(t,\xi),
\end{equation}
this means we require $F(t,\xi)-F_k(t,\xi)\in S_N\{-k-1,k+2\}$ for all $k\in\N$,
and similarly
\begin{equation}
  N(t,\xi)\sim\sum_{k=0}^\infty N^{(k)}(t,\xi),
\end{equation}
which can be chosen to be invertible,
equation \eqref{eq:II:perf_diag} implies
\begin{equation}
  \big(\D_t-\mathcal D(\xi)-R(t)\big)N(t,\xi)-N(t,\xi)\big(\D_t-\mathcal D(\xi)-F(t,\xi)\big)
  \in \bigcap_{k\in\N} S_N\{-k,k+1\}.
\end{equation}
Thus if we define the residual symbol classes
\begin{equation}
  \mathcal H\{m\}:=\bigcap_{m_1+m_2=m}S_N\{m_1,m_2\}
\end{equation}
we can find $P_\infty(t,\xi)\in\mathcal H\{1\}$ such that
\begin{equation}
  \big(\D_t-\mathcal D(\xi)-R(t)\big)N(t,\xi)=
  N(t,\xi)\big(\D_t-\mathcal D(\xi)-F(t,\xi)-P_\infty(t,\xi)\big).
\end{equation}

The classes $\mathcal H\{m\}$ are invariant under multiplication by $\exp(\pm it|\xi|)$. 
This explains why we perform more than one diagonalization step. Multiplication by
$e^{\pm it|\xi|}$ is not a well defined operation on the symbol classes 
$S_N\{m_1,m_2\}$, it destroys the symbol estimates according to the following 
proposition.

\begin{prop}\label{prop:II:H-ideal}\ \\\vspace{-3ex}
\begin{enumerate}
\item
  $e^{\pm it|\xi|} S_N^{\ell_1,\ell_2}\{m_1,m_2\}\into S_N^{\ell_1,\ell_2}\{m_1+\ell,m_2-\ell\}$ with $\ell=\ell_1+\ell_2$ 
\item
  $e^{\pm it|\xi|}\mathcal H\{m\}\into \mathcal H\{m\}$
\end{enumerate}
\end{prop}

\begin{proof}
 It suffices to prove the first statement. It holds for $a\in S_N\{m_1,m_2\}$
 \begin{align*}
  \D_t^k\D_{|\xi|}^\alpha e^{it|\xi|}a(t,\xi)&
  =\sum_{k_1+k_2=k}\sum_{\alpha_1+\alpha_2=\alpha} C_{k_1,k_2,\alpha_1,\alpha_2}
   |\xi|^{k_1}t^{\alpha_1} e^{it|\xi|} \D_t^{k_2}\D_{|\xi|}^{\alpha_2} a(t,\xi)\\
 &\le \sum_{k_1+k_2=k}\sum_{\alpha_1+\alpha_2=\alpha} C_{k_1,k_2,\alpha_1,\alpha_2}'
   |\xi|^{m_1-\alpha_2+k_1} \left(\frac1{1+t}\right)^{m_2+k_2-\alpha_1}\\
 & \le C_{k,\alpha} |\xi|^{m_1+\ell-\alpha} \left(\frac1{1+t}\right)^{m_2-\ell+k}
 \end{align*}
 for $k\le\ell_1$, $\alpha\le\ell_2$  using Leibniz rule and the definition of the
 hyperbolic zone.
\end{proof} 

\paragraph{Fundamental solution of the diagonalized system.} After performing
several diagonalization steps we want to construct the fundamental solution of
the diagonalized system
\begin{equation}\label{eq:II:CP_diag_k}
  \big(\D_t-\mathcal D(\xi)-F_{k-1}(t,\xi)-R_k(t,\xi)\big)\mathcal E_k(t,s,\xi)=0,
  \qquad\mathcal E_k(s,s,\xi)=I\in\C^{2\times2}
\end{equation}
and obtain structural properties of it.
The construction goes along the following steps:
\begin{itemize}
\item the fundamental solution $\mathcal E_0(t,s,\xi)$ to $\D_t-\mathcal D(\xi)$, 
\item influence of the main term $F^{(0)}(t,\xi)$ of $F_{k-1}(t,\xi)$,
\item influence of $F_k(t,\xi)-F^{(0)}(t,\xi)$ and $R_k(t,\xi)$.
\end{itemize}
The fundamental solution $\mathcal E_0(t,s,\xi)$ describes a 'phase function' of a
Fourier integral operator, i.e. the oscillatory behaviour of the solution multiplier. 
The main term $F^{(0)}(t,\xi)$ describes the energy decay.
Together with the other terms it constitutes a Fourier multiplier which behaves
as symbol with restricted smoothness. The number $k$ of diagonalization steps stands in 
direct connection to the smoothness properties of this symbol.

\paragraph{Step 1.} Let
\begin{equation}
  \mathcal E_0(t,s,\xi)=\exp\left\{i(t-s)\mathcal D(\xi)\right\}
  =\begin{pmatrix}e^{i(t-s)|\xi|}&\\&e^{-i(t-s)|\xi|} \end{pmatrix}
\end{equation}
such that for 
$\widetilde{\mathcal E_0}(t,s,\xi)=\frac{\lambda(s)}{\lambda(t)}\mathcal E_0(t,s,\xi)$
the equation 
\begin{equation}
  \D_t \widetilde{\mathcal E_0}(t,s,\xi)
  =\big(\mathcal D(\xi)+F^{(0)}(t,\xi)\big)\widetilde{\mathcal E_0}(t,s,\xi)
\end{equation}
is satisfied. Thus $\widetilde{\mathcal E_0}$ describes the influence of the main diagonal
term.

\paragraph{Step 2.} By the aid of $\widetilde{\mathcal E_0}(t,s,\xi)$ we define
\begin{align}
  \mathcal R_k(t,s,\xi)&=\widetilde{\mathcal E_0}(s,t,\xi)
   \big(F_{k-1}(t,\xi)+R_k(t,\xi)-F^{(0)}(t,\xi)\big)\widetilde{\mathcal E_0}(t,s,\xi),\notag\\
 &= F_{k-1}(t,\xi)+\mathcal E_0(s,t,\xi)R_k(t,\xi)\mathcal E_0(t,s,\xi)-F^{(0)}(t,\xi),
\end{align}
such that using the solution to
\begin{equation}\label{eq:CP_Q}
  \D_t\mathcal Q_k(t,s,\xi)=\mathcal R_k(t,s,\xi)\mathcal Q_k(t,s,\xi),\qquad 
  \mathcal Q_k(s,s,\xi)=I\in\C^{2\times2}
\end{equation}
the matrix $\mathcal E_k(t,s,\xi)$ can be represented as
\begin{equation}
  \mathcal E_k(t,s,\xi)=\widetilde{\mathcal E_0}(t,s,\xi)\mathcal Q_k(t,s,\xi)
  =\frac{\lambda(s)}{\lambda(t)}\mathcal E_0(t,s,\xi)\mathcal Q_k(t,s,\xi).
\end{equation}

The solution to \eqref{eq:CP_Q} is given by the Peano-Baker formula as
\begin{equation}\label{eq:Q_def}
  \mathcal Q_k(t,s,\xi)=I+\sum_{\ell=1}^\infty i^\ell
  \int_s^t \mathcal R_k(t_1,s,\xi)\int_s^{t_1}\mathcal R_k(t_2,s,\xi)\dots
  \int_s^{t_{\ell-1}}\mathcal R_k(t_\ell,s,\xi)\d t_\ell\dots\d t_1.
\end{equation}

To estimate this representation we use the following well-known inequality. It holds
\begin{prop}\label{prop:IntKerEst}
   Assume $r\in L^1_{loc}(\R)$. Then
  \begin{equation}\label{eq:IntKerEst}
    \left|\int_s^t r(t_1) \int_s^{t_1}r(t_2)\dots\int_s^{t_{k-1}}r(t_k)\d t_k\dots\d t_1 \right|
     \leq \frac1{k!}\left(\int_s^t|r(\tau)|\d\tau\right)^k
  \end{equation}
  for all $k\in\N$.
\end{prop}

\paragraph{Step 3.} The series representation \eqref{eq:Q_def} for $\mathcal Q_k(t,s,\xi)$
can be used to deduce estimates. From the unitarity of $\mathcal E_0(t,s,\xi)$ it follows
that
$$ ||\mathcal R_k(t,s,\xi)||=||R_k(t,\xi)||\lesssim \frac1{(1+t)^2|\xi|} $$
and thus using
$$ \int_{t_\xi}^\infty \frac{\d\tau}{(1+\tau)^2|\xi|}=\frac1{(1+t_\xi)|\xi|}=\frac1N $$
from the representation of $\mathcal Q_k$ by Peano-Baker formula it follows
$$
  ||\mathcal Q_k(t,s,\xi)||\lesssim 1.
$$
In a second step we want to estimate $\xi$-derivatives of 
$\mathcal Q_k(t,s,\xi)$. Proposition \ref{prop:II:H-ideal} yields for
$R_k(t,\xi)\in S_N^{\ell-k,\infty}\{-k,k+1\}$ under Assumption (A2)$_\ell$ with
$$ k-1\le\ell-k $$
that $\mathcal R_k(t,s,\xi)\in S_N^{k-1,k-1}\{-1,2\}$
uniform in $s$ (and derivations w.r.to $s$ behave like multiplications by $|\xi|$).
 
\begin{prop}\label{prop:II:peano_baker}
  Assume $a\in S_N^{\ell,\ell}\{-1,2\}$. Then
  $$
    b(t,s,\xi)=1+ \sum_{j=1}^\infty \int_s^t a(t_1,\xi)\int_s^{t_1} a(t_2,\xi)\dots
    \int_s^{t_{j-1}}a(t_j,\xi)\d t_j\dots \d t_1
  $$
  defines a symbol from $S^{\ell,\ell}_N\{0,0\}$ uniform in $s\ge t_\xi$.
\end{prop}
\begin{proof}
  We use Proposition \ref{prop:IntKerEst} to estimate this series. This yields
  in a first step (without taking derivatives)
  \begin{multline*}
    |b(t,s,\xi)|\lesssim 1+\sum_{j=1}^\infty \int_s^t \frac1{|\xi|(1+t_1)^2}
   \int_s^{t_1}\frac1{|\xi|(1+t_2)^2}\cdots\int_s^{t_{j-1}} \frac1{|\xi|(1+t_j)^2}
   \d t_j\cdots \d t_1\\
   \lesssim \exp\left\{\int_{t_\xi}^t \frac{\d\tau}{|\xi|(1+\tau)^2}\right\}\lesssim 1
  \end{multline*}
  and taking $\alpha$ derivatives with respect to $\xi$ yield in each summand 
  further factors $|\xi|^{-|\alpha|}$ according to Leibniz rule.
\end{proof}

An almost immediate consequence of this proposition is the following structural
representation of the fundamental solution.
\begin{thm}\label{thm:II:Zhyp_highreg_solrep}
  Assume (A1) and (A2)$_{2k-1}$, $k\ge1$.
  Then the fundamental solution $\mathcal E_k(t,s,\xi)$ of the diagonalized 
  system \eqref{eq:II:CP_diag_k} can be represented as
  $$ \mathcal E_k(t,s,\xi)
     =\frac{\lambda(s)}{\lambda(t)}\mathcal E_0(t,s,\xi)\mathcal Q_k(t,s,\xi) 
  \qquad t,s\ge t_\xi$$
  with a symbol $\mathcal Q_k(t,s,\xi)$ of restricted smoothness subject to the 
  symbol estimates
  $$
   \left|\left|\D_s^{\ell}\D_\xi^\alpha \mathcal Q_k(t,s,\xi)\right|\right|
   \le C_{\ell,\alpha} |\xi|^{\ell-|\alpha|}
   \qquad\qquad  t,s\ge t_\xi
  $$
  for all multi-indices $|\alpha|\le k-1$ and all $\ell\in\N_0$.
\end{thm}

The first derivative with respect to $t$ can be estimated by the equation for $\mathcal Q_k$
directly.

Of special interest is $\mathcal E_k(t,t_\xi,\xi)$. The estimate of the previous
lemma together with the properties of the derivatives of $t_\xi$ imply
\begin{cor}
  Assume (A1) and (A2)$_{2k-1}$, $k\ge1$. Then
  $$ \mathcal Q_k(t,t_\xi,\xi)\in S_N^{1,k-1}\{0,0\} $$
  (for $t\ge t_\xi$ and $|\xi|\le N$)\footnote{We could also make $t_\xi$ a smooth function in $\xi$
  to avoid such problems. But then we can not think about fundamental solutions
  $\mathcal E(t,s,\xi)$ with variable starting time $s$.}
\end{cor}

The matrix $\mathcal Q_k(t,s,\xi)$ converges for 
$t\to\infty$ to a well-defined limit. This limit will be used in Section \ref{sec:II:sharpness} to conclude the
sharpness of our results.

\begin{thm}\label{thm:II:Qinfty}
  Assume (A1) and (A2)$_{2k-1}$, $k\ge1$. The limit 
  $$
    \mathcal Q_k(\infty,s,\xi)=\lim_{t\to\infty} \mathcal Q_k(t,s,\xi)
  $$
  exists uniform in $\xi$ for $|\xi|>\xi_s$. Furthermore 
  $$
    ||\D_\xi^\alpha \mathcal Q_k(\infty,t_\xi,\xi)||\le C_\alpha |\xi|^{-|\alpha|}
  $$
  for all multi-indices $|\alpha|\le k-1$.
\end{thm}
\begin{proof}
We fix the starting value $s$ and consider only $|\xi|\ge \xi_s$ 
(i.e. $s\ge t_\xi$). Taking the difference 
$\mathcal Q_k(t,s,\xi)-\mathcal Q_k(\tilde t,s,\xi)$ in the series 
representation \eqref{eq:Q_def} yields
\begin{align*}
 \mathcal Q_k(t,s,\xi)-\mathcal Q_k(\tilde t,s,\xi)
 &=\sum_{j=1}^\infty \int_{\tilde t}^t 
  \mathcal R_k(t_1,s,\xi)\int_s^{t_1}\mathcal R_k(t_2,s,\xi)\dots
  \int_s^{t_{\ell-1}}\mathcal R_k(t_\ell,s,\xi)\d t_\ell\dots\d t_1  
\end{align*}
such that with Proposition \ref{prop:IntKerEst}
\begin{align*}
 ||\mathcal Q_k(t,s,\xi)-\mathcal Q_k(\tilde t,s,\xi)||_{L^\infty
 \{|\xi|\ge \xi_s\}}
 &\le \int_{\tilde t}^t ||R(t_1,\xi)||
   \exp\left\{\int_{s_\xi}^\infty ||R(\tau,\xi)||\d\tau\right\}\d t_1\\
 &\to 0,\qquad t,\tilde t\to\infty.
\end{align*}
Similarly one obtains for $|\alpha|\le k-1$
$$  ||\D_\xi^\alpha \mathcal Q_k(t,s,\xi)-
    \D_\xi^\alpha \mathcal Q_k(\tilde t,s,\xi)||
    \lesssim |\xi|^{-\alpha} \int_{\tilde t}^t \frac{\d\tau}{|\xi|(1+\tau)^2}
    \to0\qquad t,\tilde t\to\infty
$$
uniform in $|\xi|\ge\xi_s$
Now the second statement follows from the estimates of $t_\xi$, formula 
\eqref{eq:txi_est}.
\end{proof}

We have even proved more. The limit exists in the symbol class $\dot S^0$ of
restricted smoothness $k-1$. 

Proposition \ref{prop:IntKerEst} may also be used to estimate the formal 
representation of $\mathcal Q_k(\infty,s,\xi)$ as symbol in $(s,\xi)$.

\begin{cor}
  The series representation
  $$ 
     \mathcal Q_k(\infty,s,\xi)=
     I+\sum_{j=1}^\infty i^j \int_s^\infty  
  \mathcal R_k(t_1,s,\xi)\int_s^{t_1}\mathcal R_k(t_2,s,\xi)\dots
  \int_s^{t_{\ell-1}}\mathcal R_k(t_\ell,s,\xi)\d t_\ell\dots\d t_1
  $$
  gives an asymptotic expansion of $\mathcal Q_k(\infty,s,\xi)$
  in $S_N^{0,k-1}\{0,0\}$, i.e. the $j$th term belongs to 
  $S_N^{0,k-1}\{-j,j\}$.
\end{cor}

\paragraph{Step 4.} The transpose
of the inverse of $\mathcal Q_k$ satisfies the related equation
\begin{equation}
  \D_t\mathcal Q_k^{-T}(t,s,\xi)+\mathcal R_k^T(t,s,\xi)\mathcal Q_k^{-T}(t,s,\xi)=0,
  \qquad \mathcal Q_k^{-T}(s,s,\xi)=I\in\mathbb C^{2\times 2}.
\end{equation}
The matrix $\mathcal R_k^T(t,s,\xi)$ satisfies the same estimates like 
$\mathcal R_k(t,s,\xi)$ and therefore the reasoning of the previous step
holds in the  same way for $\mathcal Q_k^{-T}(t,s,\xi)$. Especially the
matrix $\mathcal Q_k(t,s,\xi)$ is invertible in the hyperbolic zone and
$\mathcal Q_k^{-1}(\infty,s,\xi)$ exists. 

\begin{cor}
 Assume (A1) and (A2)$_{2k-1}$, $k\ge1$. Then the limit 
  $$
    \mathcal Q_k^{-1}(\infty,s,\xi)=\lim_{t\to\infty} \mathcal Q_k^{-1}(t,s,\xi)
  $$
  exists uniform in $\xi$ for $|\xi|>\xi_s$.
\end{cor}

\paragraph{Transforming back to the original problem.} After constructing
the fundamental solution $\mathcal E_k(t,s,\xi)$ we transform back to the 
original problem and get in the hyperbolic zone the representation
\begin{equation}
  \mathcal E(t,s,\xi)=MN_k(t,\xi)\mathcal E_k(t,s,\xi) N_k^{-1}(s,\xi)M^{-1}
\end{equation}
with uniformly bounded coefficient matrices $N_k, N_k^{-1}\in S_N\{0,0\}$. We 
combine this representation with the representation obtained in the dissipative zone. This yields for $s<t_\xi$ the representation
\begin{equation}
  \mathcal E(t,s,\xi)=\frac1{\lambda(t)}
    MN_k(t,\xi)\mathcal E_0(t,t_\xi,\xi)\mathcal Q(t,t_\xi,\xi) 
    N_k^{-1}(t_\xi,\xi)M^{-1}\lambda(t_\xi)\mathcal E(t_\xi,s,\xi).
\end{equation}

\section{Estimates}

The so far obtained representations of solutions allow us to conclude estimates for the asymptotic behaviour. This section is devoted to the study of estimates
which are directly related to our micro-energy \eqref{eq:U_def}, 
i.e. estimates for the fundamental solution $\mathcal E(t,s,\D)$ 
or to the closely related energy operator $\mathbb E(t,\D)$.

\subsection{$L^2$--$L^2$ estimates}
Theorem \ref{thm:II:Zhyp_highreg_solrep} together with the estimate of Lemma \ref{lem:Zdiss1} implies
by Plancherel's theorem the following operator norm estimate.

\begin{thm} 
  Assume (A1), (A2)$_1$ and (A3). Then the $L^2$--$L^2$ estimate
  $$ ||\mathcal E(t,s,\D)||_{2\to2}\lesssim \frac{\lambda(s)}{\lambda(t)} $$
  holds. 
\end{thm}

Using the definition of the micro-energy \eqref{eq:U_def} we can reformulate 
this estimate in terms of the energy operator $\mathbb E(t,\D)$. For 
convenience we recall the relation between the multiplier $\mathcal E(t,s,\xi)$
 and $\mathbb E(t,\xi)$. They are a direct consequence of the definition of
our micro-energy \eqref{eq:U_def}. 

\begin{prop}
\begin{enumerate}
\item It holds $\mathbb E(t,\xi)=\mathcal E(t,s,\xi)\mathbb E(s,\xi)$ 
  for $s\ge t_\xi$.
\item $$ \mathbb E(t,\xi)\begin{pmatrix} \frac{h(0,\xi)}{\<\xi\>} &\\
  &1\end{pmatrix}=\begin{pmatrix} \frac{|\xi|}{h(t,\xi)} &\\
  &1\end{pmatrix}\mathcal E(t,0,\xi)$$
\item The multiplier $|\xi|/h(t,\xi)$ induces a uniformly bounded
  family of operators on $L^p$, $p\in(1,\infty)$ 
  converging strongly to the identity for
  $t\to\infty$.
\end{enumerate}
\end{prop}

\begin{cor}\label{cor:II:EnEst}
  Assumptions (A1) (A2)$_1$, (A3) imply
  $$ ||\mathbb E(t,\D)||_{2\to 2}\lesssim\frac1{\lambda(t)}. $$
\end{cor}

We conclude this section with examples. 
\begin{expl}
  Let 
  $$ b(t)=\frac\mu{1+t},\qquad \mu\in(0,1).$$
  Then Assumptions (A1) (A2) and (A3) are satisfied and the above Corollary 
  gives again the known estimate $||\mathbb E(t,\D)||_{2\to2}^2\lesssim (1+t)^{-\mu}$
  from \cite{Wir02}. 
\end{expl}
\begin{expl} \label{expl:3.5}
   Let $\mu>0$ and $n\ge 1$. Then we consider 
   $$ b(t)=\frac{\mu}{(1+t)\log(e+t)\cdots\log^{[n]}(e^{[n]}+t)}. $$
   Again the assumptions are satisfied and we obtain 
   $$ \lambda(t)=\big(\log^{[n]}(e^{[n]}+t)\big)^{\frac\mu2} $$
   and the energy decay rate may become arbitrary small. This example
   is taken from the paper of K.~Mochizuki and H.~Nakazawa,
   \cite{MN96}.
\end{expl}
\begin{expl}
  To give at least one example with an oscillating coefficient we consider
  $$
    b(t)=\frac{2+\cos(\alpha\log(e+t))}{4(e+t)}
  $$
  with large real $\alpha$. Then it holds
  $$ \int b(t)\d t=\frac12\log(e+t)+\frac1{4\alpha}\sin(\alpha\log(e+t)) $$
  and thus
  $$ \lambda(t)=C_0\sqrt[4]{e+t}\exp\left\{\frac1{4\alpha}\sin(\alpha\log(e+t))\right\} $$
  such that the energy decay is given by $||\mathbb E(t,\D)||_{2\to2}^2\sim (e+t)^{-\frac12}$,
  which is independent on the choice of $\alpha$.
\end{expl}

\subsection{$L^p$--$L^q$ estimates}

This section is devoted to $L^p$--$L^q$ estimates.  The basic estimate
is given in the following theorem, it restates the known result for free
waves in the language of our operators. 
\begin{thm} It holds
  $$ ||\mathcal E_0(t,\D)||_{p,r\to q}\le C_{p,q} (1+t)^{-\frac{n-1}2\left(\frac1p-\frac1q\right)} $$
  for dual indices $p$ and $q$, $p\in(1,2]$ and with regularity $r=n\left(\frac1p-\frac1q\right)$.
\end{thm}

By the aid of this estimate we deduce from our representation a corresponding
estimate for the dissipative equation.

\begin{thm}\label{thm:II:LpLq-decay}
  Assume (A1), (A2) and (A3). Then the operator $\mathcal E(t,s,\D)$ satisfies
  for dual indices $p$ and $q$, $p\in(1,2]$, $pq=p+q$ the norm estimate
  $$
    ||\mathcal E(t,0,\D)||_{p,r\to q}\lesssim
     \frac1{\lambda(t)} (1+t)^{-\frac{n-1}2\left(\frac1p-\frac1q\right)}
  $$
  with regularity $r=n\left(\frac1p-\frac1q\right)$.
\end{thm}
\begin{proof}
We decompose the proof into two parts. First we consider $\mathcal E(t,0,\D)\phi_{diss,N}(t,\D)$. 
Using the estimate $||\mathcal E(t,0,\xi)\phi_{diss,N}(t,\xi)||\lesssim \frac1{\lambda^2(t)}$ together
with the definition of the zone we get
$$
  ||\mathcal E(t,0,\D)\phi_{diss,N}(t,\D)||_{p\to q}
     \lesssim \frac1{\lambda^2(t)} (1+t)^{-n\left(\frac1p-\frac1q\right)}
$$
which is a stronger decay rate than the one given in the theorem.

In a second step we consider the hyperbolic part. For small frequencies we use the representation 
\begin{multline*}
    \mathcal E(t,0,\xi)\phi_{hyp,N}(t,\xi)=\\
   \frac1{\lambda(t)}\underbrace{MN_k(t,\xi)}_{q\to q}
   \underbrace{\mathcal E_0(t-t_\xi,\xi)}_{p,r\to q}
   \underbrace{\mathcal Q_k(t,t_\xi,\xi)}_{p,r\to p,r}
   \underbrace{N_k^{-1}(t_\xi,\xi)M^{-1}\lambda(t_\xi)\mathcal E(t_\xi,0,\xi)}_{p,r\to p,r}
   \phi_{hyp,N}(t,\xi)
\end{multline*}
together with the mapping properties of the multipliers marked with a brace. It is essential that
$k-1\ge \lceil\frac n2\rceil$. The operator $\mathcal E_0(t,\xi)$ brings the hyperbolic decay rate, while
$\mathcal E_0(-t_\xi,\xi)\in S^0$.

For large frequencies the representation simplifies to
$$
    \mathcal E(t,0,\xi)\phi_{hyp,N}(t,\xi)=
   \frac1{\lambda(t)}\underbrace{MN_k(t,\xi)}_{q\to q}
   \underbrace{\mathcal E_0(t,\xi)}_{p,r\to q}
   \underbrace{\mathcal Q_k(t,0,\xi)}_{p,r\to p,r}
   \underbrace{N_k^{-1}(0,\xi)M^{-1}}_{p,r\to p,r}
   \phi_{hyp,N}(t,\xi),
$$
the argumentation remains the same.
\end{proof}

\paragraph{Minimal regularity for the $L^p$--$L^q$ estimate.} With the notation
$$
   \ell_n=2\lceil\frac n2\rceil+1=\begin{cases}n+1,\qquad& \text{$n$ even},\\n+1,&\text{$n$ odd,}\end{cases}
$$
we can prove the above  $L^p$--$L^q$ decay estimate under the weaker Assumption 
(A2)$_{\ell_n}$ on the coefficient function. If we use this regularity of the coefficient
and perform $k=\lceil\frac n2\rceil$ diagonalization steps we obtain
$ N_k(t,\xi)\in S_N^{0,\infty}\{0,0\}$
and $\mathcal Q_k(t,s,\xi)$ is uniformly in $t\ge s\ge t_\xi$ a symbol of smoothness 
$\lceil\frac n2\rceil$. Thus, $N_k\phi_{hyp}$, $N_k^{-1}\phi_{hyp}$ and 
$\mathcal Q_k\phi_{hyp}$ define operators
$L^p\to L^p$ for all $p\in(1,\infty)$ with uniformly bounded operator norm in $t, s$.

\subsection{Estimates of the solution itself}

The definition of the micro-energy allows to extract also estimates for the solution from
the structure of $\mathcal E(t,s,\xi)$. Similar to the reasoning in the previous theorems
one obtains for the solution operator
$$
  \mathbb S(t,\D)\;:\;(u_1,\<D\>^{-1} u_2) \mapsto u(t,\cdot),
$$
normalised in such a way that $\mathbb S(t,\D):L^2\to L^2$,
the following $L^p$--$L^q$ estimate.

\begin{thm}\begin{enumerate}\item
  Assume (A1), (A2)$_1$ and (A3). Then the solution $u(t,x)$ satisfies the 
  $L^2$--$L^2$ estimate
  $$
     ||\mathbb S(t,\D)||_{2\to2} \lesssim \frac{1+t}{\lambda^2(t)}.
  $$
  \item 
  Assume (A1), (A2)$_{\ell_n}$ and (A3). Then the solution $u(t,x)$ satisfies the 
  $L^p$--$L^q$ estimate
  $$
     ||\mathbb S(t,\D)||_{p,r\to q} \lesssim \begin{cases}
    \frac1{\lambda^2(t)}(1+t)^{1-n\left(\frac1p-\frac1q\right)},\qquad &p\ge p^*,\\ 
    \frac1{\lambda(t)} (1+t)^{-\frac{n-1}2\left(\frac1p-\frac1q\right)}, &p< p^*,
   \end{cases}
  $$
  for dual indices $p\in(1,2]$, $pq=p+q$ and with regularity $r=n\left(\frac1p-\frac1q\right)$. 
  The critical value $p^*=\frac{n+3}{2n+2}$ if $tb(t)\to0$ and determined by the balance of the two estimates 
  in the remaining cases. 
 \end{enumerate}
\end{thm}
\begin{proof} 
  We start with the $L^2$--$L^2$ estimate.  

  In the hyperbolic zone and for $|\xi|>c$ we estimated $|\xi|\hat u$ by $\lambda^{-1}(t)$ and we 
  can just divide by $|\xi|$ to get an estimate for $u$ by 
  $$\frac1{\lambda(t)}\lesssim\frac{1+t}{\lambda^2(t)}. $$ 

  For small $|\xi|$ we have to take into account that dividing by $h(t,\xi)$ 
  brings a further factor $(1+t)$. Thus we obtain inside the dissipative zone $(1+t)/\lambda^2(t)$
  and in the hyperbolic zone for small frequencies
  $$
    \frac1{|\xi|\lambda(t)\lambda(t_\xi)}\sim \frac{1+t_\xi}{\lambda(t)\lambda(t_\xi)}
    \lesssim \frac{1+t}{\lambda^2(t)}
  $$
  using the monotonicity of $t/\lambda(t)$ for large $t$ following from (A3). 

  To conclude the $L^p$--$L^q$ estimate we follow the same proof like in the case for the energy
  with the exception that the rate inside the dissipative zone is different and dominating
  at least if $p$ and $q$ are near to two.
\end{proof}

This result coincides for the case $b(t)=\frac{\mu}{1+t}$, $\mu<1$ with the 
estimate from \cite{Wir02}. We think the following observation is important
to understand the essential difference between the solution estimate and the estimate
for the energy. 

\begin{quote}
  For $p$ near to $2$ the estimate for the solution operator
  $\mathbb S(t,\D)$ comes from properties of the dissipative zone, while estimates for
  the energy $\mathbb E(t,\D)$ are determined by the hyperbolic zone.
\end{quote}

\begin{rem}
  Assumptions on the data, which are effective near the frequency $\xi=0$ may be used to
  improve these estimates. One possibility to achieve this is the use of weighted initial
  data. If we assume that $\<x\>u_1,\<x\>u_2\in L^2$ and the space dimension $n\ge 3$, then
  the estimate inside the dissipative zone in this theorem might be improved by an application
  of Hardy inequality and one obtains
  $$ ||\mathbb S(t,\D)||_{\<x\>^{-1}L^2\to L^2}\lesssim \frac1{\lambda(t)}. $$
\end{rem}

\section{Sharpness}\label{sec:II:sharpness}
Finally we want to prove the sharpness of the above given energy decay 
estimates.
Our constructive approach enables us to formulate the question of sharpness
as a \eIndex{modified scattering theory}. The basic idea is as follows
\begin{itemize}
\item
  we relate the energy operator $\mathbb E(t,D)$ to the
  corresponding operator $\mathbb E_0(t,\D)$ for free waves
  multiplied by the decay rate 
\item
  this relation defines a M\o{}ller wave operator $W_+$ defining
  appropriate data to the free wave equation with the same asymptotic
  properties (up to the factor $\lambda(t)$),
\item
  furthermore, we need to know the mapping properties of the 
  M\o{}ller wave operator,
\item
  and the convergence defining the wave operator has to 
  be understood.
\end{itemize}

A first observation follows immediately from Liouville Theorem
 and gives an expression for the determinant of 
$\mathbb E(t,\xi)$.

\begin{lem}\label{lem:II:det_bbE}
  It holds $\det\mathbb E(t,\xi)=\frac1{\lambda^2(t)}[\xi]$, where $[\xi]=|\xi|/\<\xi\>$.
\end{lem}

After these introductory remarks we can state the following theorem. It holds
\begin{thm}\label{thm:II:mod_scattering}
  Assume (A1), (A2) and (A3). Then the limit
  $$ W_+(\D)=\slim_{t\to\infty} \lambda(t)(\mathbb E_0(t,\D))^{-1}
    \mathbb E(t,\D) $$
  exists as strong limit in $L^2\to L^2$ and defines the \eIndex{modified
  M\o{}ller wave operator} $W_+$. It satisfies
  $$
   W_+(\xi)=(\mathbb E_0(t_\xi,\xi))^{-1}
   \mathcal Q_k(\infty,t_\xi,\xi)N_k^{-1}(t_\xi,\xi)M^{-1}
   \lambda(t_\xi) \mathbb E(t_\xi,\xi)
  $$ 
  for all $k\ge 1$.
\end{thm}

Remark that $t_\xi$ depends on the zone constant and this constant
is chosen after diagonalizing $k$ steps. Thus 
$$
  \mathcal Q_k(\infty,t_\xi,\xi)N_k^{-1}(t_\xi,\xi)
$$
is independent on $1\le k\le\ell$ for sufficiently large zone constant $N$
depending on $\ell$. 

\begin{proof}
The proof consists of three steps.
\paragraph{Step 1.} With the notation
$$
   V_c=\{U\in L^2|\dist(0,\supp \hat U)\ge c\}
$$
we can construct the dense subspace $M=\bigcup_{c>0} V_c$. Now Theorem
\ref{thm:II:Qinfty} together with the representation 
$\mathbb E(t,\xi)=\mathcal E(t,t_\xi,\xi)\mathbb E(t_\xi,\xi)$
implies the existence of the limit 
$$ 
  \lim_{t\to\infty}\lambda(t) \mathbb E_0(t,\D)^{-1}\mathbb E(t,\D)
$$
as limit in the operator norm in $V_c\to V_c$ for all $c>0$. Thus the limit exists pointwise
on $M$.

\paragraph{Step 2.} The energy estimate, Theorem \ref{cor:II:EnEst},
implies that $\lambda(t)\mathbb E_0(t,\D)^{-1}\mathbb E(t,\D)$ is uniformly
bounded in $L^2\to L^2$. Thus the Theorem of Banach-Steinhaus implies the 
existence of the strong limit and defines $W_+$.

\paragraph{Step 3.} The previously defined operator $W_+$ is given on each
subspace $V_c$ as Fourier multiplier with symbol
$$
  W_+(\xi)=(\mathbb E_0(t_\xi,\xi))^{-1}
   \mathcal Q_k(\infty,t_\xi,\xi)N_k^{-1}(t_\xi,\xi)M^{-1}
   \lambda(t_\xi) \mathbb E(t_\xi,\xi),
$$
which is independent of $c$. Thus the representation holds on $M$ and 
using the boundedness of $W_+$ on the whole space.
\end{proof}

Remark that for above given theorem we need (A2)$_{2k-1}$ to represent $W_+(\xi)$ in terms of
$\mathcal Q_{k}(t,s,\xi)$. Thus Assumption (A2)$_1$ is sufficient to get the existence of the
wave operator.

\begin{cor}\label{cor:II:scatt1}
  It holds
  $$
     \det W_+(\xi)=1,
  $$
  and therefore $W_+(\D)$ defines an isomorphism on $L^2$.
\end{cor}

\paragraph{Interpretation of the result.}  Theorem
\ref{thm:II:mod_scattering} may be used to construct for each datum 
$(\<\D\>u_1,u_2)\in L^2$ to Cauchy problem \eqref{chap3:eq:CP} 
a corresponding datum $(\<\D\>\tilde u_1,\tilde u_2)^T=W_+(\D)(\<\D\>u_1,u_2)^T$
to the free wave equation $\bx\tilde u=0$, such that the solutions
are asymptotically equivalent up to the decay factor $\lambda^{-1}(t)$, i.e.
it holds
\begin{equation}
  \left|\left|\mathbb E_0(t,\D)(\<\D\>\tilde u_1,\tilde u_2)^T
   -\lambda(t)\mathbb E(t,\D)(\<\D\>u_1,u_2)^T\right|\right|_2
  \to0\qquad \text{as $t\to\infty$}.
\end{equation}
This is a direct consequence of the unitarity of $\mathbb E_0(t,\xi)$. It 
implies that the above given $L^2$--$L^2$ estimates are indeed sharp and 
describe for all nonzero initial data the exact decay rate.

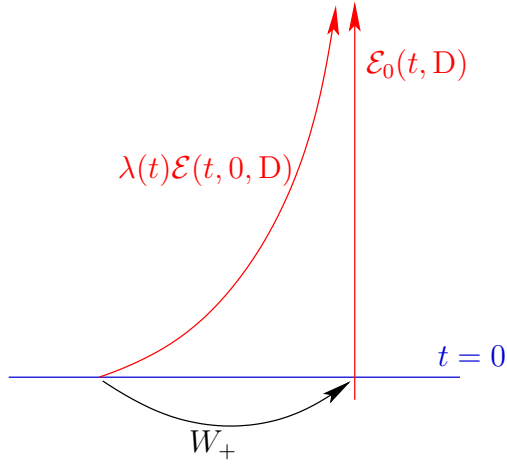
\begin{figure}
\begin{center}
\input{fig_scatt_mod.pstex_t}
\caption{Modified scattering theory.}
\end{center}
\end{figure}

\begin{appendix}
\section{Appendix: A remark on Volterra integral equations}   
We are interested in solutions to the Volterra equation
\begin{equation}\label{eq:II:volt_int_eq}
  f(t,p)=f(0,p)+\int_0^t k(t,\tau,p) f(\tau,p)\d\tau
\end{equation}
with kernel $k=k(t,\tau,p)$ and given initial value $f(0,p)$ depending on some parameter 
$p\in P\subseteq \R^n$.

\begin{thm}\label{thm:B4}
  Assume $f(0,\cdot)\in L^\infty(P)$, $k\in L^\infty(\R_+^2\times P)$ and 
  $$ \int_0^t |k(t,\tau,p)|\d\tau\in L^\infty(\R_+\times P). $$
  Then there exists a (unique) solution $f(t,p)$ of \eqref{eq:II:volt_int_eq} in 
  $L^\infty(\R_+\times P)$.
\end{thm}

For the proof of this statement we follow \cite{GLS90}. The condition on the kernel implies that
the operator
$$ K\; :\; f(t,p)\mapsto \int_0^t k(t,\tau,p)f(\tau,p)\d\tau $$
is bounded. This implies that, if the norm of this operator is sufficiently small, the solution
of the integral equation is given by the contraction mapping principle. Now the idea of 
\cite[Chapter 9.3, Theorem 3.13]{GLS90} applies, we can decompose the time interval into a 
finite number of smaller sub-intervals where the restricted operators are contractions. 
This follows from
$$ \esssup_{p\in P} \int_{T_1}^{T_2} |k(t,\tau,p)|\d\tau
 \le (T_2-T_1)\,||k||_\infty \to 0,\qquad\qquad \text{as $T_2-T_1\to0$.} $$
Furthermore, the Volterra structure of the equation can be used to build up the resolvent
from the resolvents of the restricted operators like in \cite{GLS90}.
\end{appendix}

%\bibliographystyle{alpha}
%\bibliography{../database}

\end{document}

%% file: fig_scatt_mod.pstex_t
\begin{picture}(0,0)%
\includegraphics{fig_scatt_mod.pstex}%
\end{picture}%
\setlength{\unitlength}{2072sp}%
\begingroup\makeatletter\ifx\SetFigFont\undefined%
\gdef\SetFigFont#1#2#3#4#5{%
  \reset@font\fontsize{#1}{#2pt}%
  \fontfamily{#3}\fontseries{#4}\fontshape{#5}%
  \selectfont}%
\fi\endgroup%
\begin{picture}(6028,5530)(1779,-5029)
\put(6931,-3886){\makebox(0,0)[lb]{\smash{{\SetFigFont{12}{14.4}{\familydefault}{\mddefault}{\updefault}{\color[rgb]{0,0,.82}$t=0$}%
}}}}
\put(6076,-376){\makebox(0,0)[lb]{\smash{{\SetFigFont{12}{14.4}{\familydefault}{\mddefault}{\updefault}{\color[rgb]{1,0,0}$\mathcal E_0(t,\D)$}%
}}}}
\put(3961,-4921){\makebox(0,0)[lb]{\smash{{\SetFigFont{12}{14.4}{\familydefault}{\mddefault}{\updefault}{\color[rgb]{0,0,0}$W_+$}%
}}}}
\put(3106,-1636){\makebox(0,0)[lb]{\smash{{\SetFigFont{12}{14.4}{\familydefault}{\mddefault}{\updefault}{\color[rgb]{1,0,0}$\lambda(t)\mathcal E(t,0,\D)$}%
}}}}
\end{picture}%